\documentclass[reqno,12pt]{amsart}
\usepackage{amscd}
\usepackage{amssymb}
\usepackage[a4paper]{geometry}
\usepackage{amstext}
\usepackage{mathrsfs}
\usepackage{todonotes}
\usepackage[final]{hyperref}
\hypersetup{unicode= false, colorlinks=true, linkcolor=blue,
anchorcolor=blue, citecolor=red, filecolor=red, menucolor=blue,
pagecolor=blue, urlcolor=blue} 
\numberwithin{equation}{section}
\usepackage[none]{hyphenat}

\newcommand{\half}{\frac{1}{2}}

\newcommand{\CC}{\mathbb {C}}
\newcommand{\Hg}{\mathscr{H}}

\newtheorem{thm}{Theorem}[section]

\newtheorem{lemma}{Lemma}[section]
\newtheorem{corollary}{Corollary}[section]

\title[Compact Carleson measures
 from sparse sequences] { Compact Carleson measures from sparse
 sequences}
\author [Tesfa  Mengestie]{Tesfa  Mengestie }
\address{Department of Mathematical Sciences\\
Norwegian University of Science and Technology (NTNU)\\
 NO- 7491 Trondheim, Norway}
\email{mengesti@math.ntnu.no}
\thanks{The author is  supported by the Research Council of
Norway grant 185359/V30.} \subjclass[2000]{30E05, 46E22}

\begin{document}
\begin{abstract}
In \cite{BMS2}, Y. Belov, K. Seip,  and the author studied the
Carleson measures for certain spaces of analytic functions of which
the de Branges spaces and the model subspaces of the Hardy space
$H^2$ are the prime examples. In this paper, we  continue this line
of research by studying  the compact Carleson measures for such
spaces.
\end{abstract}
 \maketitle
\section{Introduction}
This paper is concerned with the compact embedding maps  induced by
Carleson measures  in certain spaces of analytic functions. Bounded,
unitary, and invertible maps were studied first in
  \cite{BMS2,BMS}.

  Let
 $\Gamma=(\gamma_n)$ be a sequence of
distinct complex numbers and $v=(v_n)$ be a weight sequence that
satisfies the admissibility condition
\begin{equation}
\sum_{n=1}^\infty \frac{v_n}{1+|\gamma_n|^2}<\infty. \label{adm}
\end{equation} When \eqref{adm} holds, we will call $v$ an admissible weight sequence for
  $\Gamma$. Any such pair $(\Gamma, v)$ parameterizes the space $\Hg(\Gamma,v)$ which
  consists of all functions
\[ f(z)=\sum_{n=1}^\infty\frac{a_n v_n}{z-\gamma_n} \]
for which
\[ \|f\|_{\Hg(\Gamma,v)}= \|(a_n)\|_{\ell_v^2}<\infty, \ \  \ell_v^2=
\big\{ (a_n) : \|a\|_v^2=\sum_{n=1}^\infty |a_n|^2v_n
<\infty\big\}\] and $z$ belongs to the set
 \[ (\Gamma,
v)^{*}=\Big\{ z\in \CC:\ \sum_{n=1}^\infty
\frac{v_n}{|z-\gamma_n|^2}<\infty\Big\} .\]
 It means that we obtain the value of a function $f$ in
$\Hg(\Gamma, v)$ at a point $z$ in $(\Gamma, v)^{*}$ by computing
the
 weighted discrete Hilbert transform:
\begin{equation}
 (a_n)\mapsto\sum_{n=1}^\infty \frac{a_n
v_n}{z-\gamma_n}, \label{discreteH}
\end{equation}
which is well defined whenever $(a_n)$ belongs to $\ell_v^2.$
   The de Branges
 spaces, model subspaces of the Hardy space  $H^2$ which admit Riesz bases of reproducing kernels
  and the Fock-type spaces
 studied in \cite{BL} are  all examples of spaces of the kind  $\Hg(\Gamma, v).$  We refer to
\cite{BMS2, BMS, TYM} for detailed accounts.

We say that a nonnegative measure $\mu$ on $(\Gamma,v)^{*}$ is a
Carleson measure for $\Hg(\Gamma,v)$ if there is a positive constant
$C$ for which
\[\int_{(\Gamma,v)^*}|f(z)|^2d\mu(z)\le C \|f\|^2_{\Hg(\Gamma,v)}\]
holds for every $f$ in $\mathscr{H}(\Gamma,v)$. Such measures have
been described  in \cite{BMS2} (see Theorem \ref{3equivalent} below)
when $\Gamma$ grows at least exponentially, i. e., when
\begin{equation}
\inf_{n\geq1}|\gamma_{n+1}|/|\gamma_n|>1. \label{expon}
\end{equation}
To state the result in \cite{BMS2}, we  first partition the complex
plane $\CC$ in the following way. Set $\Omega_1=\{z\in \CC:\
|z|<\big(|\gamma_1|+|\gamma_2|\big)/2\}$ and then \[
\Omega_n=\left\{z\in \mathbb{C} : \
\big(|\gamma_{n-1}|+|\gamma_n|\big)/2\le
|z|<\big(|\gamma_n|+|\gamma_{n+1}|\big)/2\right\}
\text{ for}\ \ n\ge 2.\] 
\begin{thm} [Y. Belov, T. Mengestie, K. Seip \cite{BMS2} ]\label{3equivalent}
Suppose that the sequence $\Gamma$ satisfies the sparseness
condition \eqref{expon} and that $v$ is an admissible weight
sequence for $\Gamma$. If $\mu$ is a nonnegative measure on $\CC$
with $\mu(\Gamma)=0$, then the following are equivalent.
\begin{enumerate}
\setlength{\itemsep}{0ex}
\item The weighted Hilbert transform in \eqref{discreteH} is bounded from $\ell^2_v$ to
$L^2(\mathbb{C},\mu).$
\item The measure $\mu$ is a Carleson measure for the space $\Hg(\Gamma, v)$.
\item We have that
\begin{equation*}
\label{trivial} \sup_{n\geq1} \int_{\Omega_n} \frac{v_n
d\mu(z)}{|z-\gamma_n|^2}<\infty
\end{equation*} and
\begin{equation*}
\label{carleson} \sup_{n\geq1} \left( \sum_{m=1}^{n}v_m
\sum_{m=n+1}^\infty \int_{\Omega_m} \frac{d\mu(z)}{|z|^2} +\sum_{m=
n+1}^\infty \frac{v_n}{ |\gamma_n|^2}
\sum_{m=1}^{n}\mu(\Omega_m)\right)<\infty.
\end{equation*}
\end{enumerate}
\end{thm}
The purpose  of this paper is to identify  those Carleson measures
$\mu$ for which the embedding maps $ I_\mu$ from $\Hg(\Gamma, v)$
into $L^2(\CC, \mu)$ are  compact. Whenever $\mu$ induces such an
embedding, we call it a vanishing or compact Carleson measure for
$\Hg(\Gamma, v).$
Compact  Carleson measures appear naturally in the study of compact
composition operators in various function spaces. As far as  their
characterization is concerned, there  exists a  general ``folk
theorem'':
 once the Carleson measures are described by a certain ``big oh''
condition, vanishing Carleson measures are then characterized by the
corresponding ``little oh'' counterparts. This does not mean that
such  `` folk theorem'' is always true. See  \cite{FMB} for a
counterexample. Our result shows that  it  still holds in the space
$\Hg(\Gamma,v)$.
\begin{thm} \label{compthm}
Suppose that the sequence $\Gamma$ satisfies the sparseness
condition \eqref{expon} and that $v$ is an admissible weight
sequence for $\Gamma$. A nonnegative measure $\mu$ on $\CC$ with
$\mu(\Gamma)=0$
  is a compact Carleson measure for $\Hg(\Gamma,v)$ if and only if
\begin{equation}
\label{trivialcomp} \lim_{n\to \infty} \int_{\Omega_n} \frac{v_n
d\mu(z)}{|z-\gamma_n|^{2}}=0
\end{equation}
and
\begin{equation}
\label{compactcarleson}\lim_{n\to \infty}
\left(\sum_{m=1}^{n}v_m\sum_{m=n+1}^\infty \int_{\Omega_m}
\frac{d\mu(z)}{|z|^{2}} + \sum_{m= n+1}^\infty\frac{ v_m}{
|\gamma_m|^{2}}\sum_{m=1}^{n}\mu(\Omega_m)\right)= 0.
\end{equation}
\end{thm}
 Condition  \eqref{trivialcomp} of Theorem \ref{compthm} is a
condition on the local behavior of $\mu$, while condition
\eqref{compactcarleson} deals with its global behavior. Combining
the two conditions, we see that \eqref{trivialcomp} may be replaced
by a stronger global necessary condition:
 \begin{equation*}
 \label{specialcase}
\lim_{n\to \infty}\int_{\CC} \frac{v_n d\mu(z)}{ |z-\gamma_n|^2} =0.
\end{equation*} We consider
 two cases in which \eqref{compactcarleson} is
  automatically fulfilled once either this condition or the
original one \eqref{trivialcomp} holds.
\begin{corollary} \label{cor1}
Suppose the sequence $\Gamma$ satisfies the sparseness condition
\eqref{expon} and that the numbers $v_n$ grow at least exponentially
and that the numbers $v_n/|\gamma_n|^2$ decay at least exponentially
with $n$. A nonnegative measure $\mu$ on $\CC$ with $\mu(\Gamma)=0$
  is a compact Carleson measure for $\Hg(\Gamma,v)$ if and only if
 \begin{equation*}
 \label{cor11}
\lim_{n\to \infty}\int_{\Omega_n}\frac{v_n d\mu(z)}{
|z-\gamma_n|^2}=0.
\end{equation*}
\end{corollary}
\begin{corollary} \label{cor11}
Suppose the sequence $\Gamma$ satisfies the sparseness condition
\eqref{expon} and that $(v_n) \in \ell_1$. A nonnegative measure
$\mu$ on $\CC$ with $\mu(\Gamma)=0$ is a compact Carleson measure
for $\Hg(\Gamma,v)$ if and only if
  \begin{equation*}
\label{cor2} \lim_{n\to \infty}\int_{\CC}\frac{v_n d\mu(z)}{
|z-\gamma_n|^2} =0.\end{equation*}
\end{corollary}
Both corollaries follow immediately from Theorem \ref{compthm}.\\
 We note that the discrete version of  Theorem \ref{3equivalent} was
  discussed in detail in \cite{BMS2,TYM} since it  describes the Bessel sequences
  of normalized reproducing kernels  in $\Hg(\Gamma,v).$ Since a
  Bessel sequence of normalized reproducing kernels corresponds to a
  compact measure if and only if the sequence is finite, discrete
  versions of Theorem \ref{compthm} are less interesting and will
  not be considered in this paper.

 An interesting question is when a compact Carleson measure $\mu$
 induces a Schatten  $p$ class embedding map $I_\mu$ from
$\Hg(\Gamma, v)$ into $L^2(\CC, \mu)$ for $ 0< p<\infty.$   For
$p=2$ (Hilbert--Schmidt class), the answer follows from the fact
that
\begin{equation}
\label{hilbertschmidt}\|I_\mu\|_{2}^2
    =\int_{\CC}\|k_z\|_{\Hg(\Gamma,v)}^2 \ d\mu(z) <\infty, \end{equation}
where $ \|I_\mu\|_{2}$ denotes the Hilbert--Schmidt norm of $I_\mu.$
We observe that the reproducing kernel of $\Hg(\Gamma,v)$ at a point
$\lambda$ in $(\Gamma, v)^*$  is explicitly given by \[k_\lambda(z)=
\sum_{n=1}^\infty
\frac{v_n}{(\overline{\lambda}-\overline{\gamma_n})(z-\gamma_n)};\]
this is a direct consequence of the definition of $\Hg(\Gamma,v)$.
Then,  a simple computation along with \eqref{expon}
 gives that \eqref{hilbertschmidt} holds if and
only if
\begin{equation*}
\label{trivialHS} \sum_{m=1}^\infty \int_{\Omega_m} \frac{v_m
d\mu(z)}{|z-\gamma_m|^2} <\infty
\end{equation*}
and
\begin{equation*} \label{BHS}
\sum_{m=1}^\infty v_m \sum_{j=m+1}^\infty \int_{\Omega_j}
\frac{d\mu(z)}{|z|^2}+  \sum_{m=1}^\infty
\frac{v_m}{|\gamma_m|^2}\sum_{j=1}^{m-1}\mu(\Omega_j)<\infty.
\end{equation*}
 For other ranges of $p,$ a few partial
 results are available in \cite{TYM}. In general, a comprehensive solution
 is yet to be found.

We close this section with a few words on  notation. The notation
$U(z)\lesssim V(z)$ (or equivalently $V(z)\gtrsim U(z)$) means that
there is a constant $C$ such that $U(z)\leq CV(z)$ holds for all $z$
in the set in question, which may be a Hilbert space, a set of
complex numbers, or a suitable index set. We write $U(z)\simeq V(z)$
if both $U(z)\lesssim V(z)$ and $V(z)\lesssim U(z)$.
\section{Proof of Theorem \ref{compthm}}
 We note that since
$\Hg(\Gamma, v)$ is reflexive,
 $\mu$ induces a
 compact embedding if and only if  each weakly convergent
sequence in $\Hg(\Gamma, v)$ converges  in norm in $L^2(\CC, \mu)$.
 It is known that a weakly convergent sequence is uniformly norm bounded.
In general, the converse statement does not hold. But under an
additional assumption, the following particular case of Nordgren's
 result holds \cite{Nord}.
 \begin{lemma}\label{nord}
 Let $(f_n)$ be a sequence of functions in $\Hg(\Gamma, v).$ Then $(f_n)$
 converges weakly to zero (weekly null)
 if and only if it
 converges pointwise to zero  and
 \begin{equation}\label{nordeq}
 \sup_n \|f_n\|_{\Hg(\Gamma,v)} <\infty.
 \end{equation}
 \end{lemma}
\subsection{Proof of the necessity of  Theorem \ref{compthm}:} For
simplicity, we set
\begin{equation*}
V_n= \sum_{m=1}^{n-1} v_m \ \ \ \text{and} \ \ \ P_n=
\sum_{m=n+1}^\infty v_m |\gamma_m|^{-2}.
\end{equation*}
  We may first choose a sequence of normalized test functions
 \[q_n(z)=\frac{ \sqrt{v_n}}{z-\gamma_n}.\]  The  sequence
converges weakly to zero in $\Hg(\Gamma, v).$ This is a particular
case of a much more general statement which says that any
orthonormal sequence in a Hilbert space converges weakly to zero.
This along with compactness of $\mu$
  yields
  \begin{equation*}
  0= \lim_{n\to \infty} \int_{\CC} |q_n(z)|^2 d\mu(z)
  \end{equation*} from which  condition \eqref{trivialcomp} follows.

We next  consider  another
  sequence of unit norm   functions defined by
 \[g_n(z)= \frac{1}{\sqrt {P_n}} \sum_{m=n+1}^\infty \frac{v_m}{\overline{\gamma_m}
(z-\gamma_m)}. \]
 If $z$ belongs to the shell $\Omega_N,$ then
$|g_n(z)|\simeq P_n^\half$ whenever $n>N,$ and therefore $g_n$
converges pointwise to zero as $n\to \infty.$
 Thus by the above lemma, the sequence $g_n$ is  weakly null.
  Taking into account the compactness of $\mu$, we find that
 \begin{equation}
 \label{tail}
  0=\lim_{n\to \infty} \int_{\CC} |g_n(z)|^2 d\mu(z)\geq \lim_{n\to \infty}
 P_n \sum_{m=1}^{n}\mu(\Omega_m).
 \end{equation}
 On the other hand, if $\sup_n V_n <\infty,$ then
\eqref{compactcarleson} holds  trivially  for each  Carleson measure
$\mu.$ We shall thus consider  the  case when $V_n \to \infty $ as
$n \to \infty$. To this end, we consider a third  sequence of
unit norm functions
 \[h_n(z)= \frac{1}{V_n^\half} \sum_{m=1}^{n-1} \frac{v_m}{z-\gamma_m}.\]
 We shall verify that $h_n$ converges pointwise
 to zero. If  $z$ belongs to the annulus $\Omega_N,$ then the estimate
\begin{equation}
\label{hardp} |h_n(z)|\simeq  \frac{1}{V_n^\half} \sum_
{m=1}^{N-1}\frac{v_m}{|z|}
 +\frac{v_N}{V_n^\half|z-\gamma_N|}+
\frac{1}{V_n^\half} \sum_ {m=N+1}^{n-1}\frac{v_m}{|\gamma_m|}
\end{equation} holds when $n$ is sufficiently large.
 It suffices to look at only the right-hand side of \eqref{hardp} because
  the first two terms clearly converge to zero
 as $n\to \infty.$ Let $c_n$ be a sequence increasing to infinity such that
 $$\sum_{j=1}^\infty  \frac{c_jv_j}{|\gamma_j|^2} <\infty.$$
 Then an application of the  Cauchy--Schwarz inequality gives
\begin{eqnarray}\frac{1}{V_n^\half} \sum_ {m=N+1}^{n-1}\frac{v_m}{|\gamma_m|}
&\leq& \left(\sum_ {m=N+1}^{n-1}\frac{c_m v_m}{|\gamma_m|^2}\right)^\half \left(\sum_{m=1}^{n-1}
\frac{v_m}{c_m V_n}\right)^\half\nonumber\\
 &\lesssim&\frac{1}{V_n^\half}\left(\sum_{m=1}^{n-1}\frac{v_m}{c_m}\right)^\half \to 0  \end{eqnarray}
 as $n \to \infty.$
    It follows that
 by Lemma \ref{nord} $h_n$ converges weakly to zero.\\
 To this end, if $\mu$ induces a compact embedding, we
then have
 \begin{equation*}
 \label{remaining}
 0= \lim_{n\to \infty} \int_{\CC} |h_n(z)|^2 d\mu(z) \gtrsim  \lim_{n\to \infty}  V_n\sum_{k=n+1}^\infty
 \int_{\Omega_k} \frac{d\mu(z)}{|z|^2},
 \end{equation*} which, together with \eqref{tail}, gives the assertion in
 \eqref{compactcarleson}.
\subsection{Proof of the sufficiency of Theorem \ref{compthm}:} Assume conversely
that the conditions \eqref{trivialcomp}, and \eqref{compactcarleson}
 hold, and consider a weak null  sequence
\[f_n(z) = \sum_{m=1}^\infty\frac{a_m^n v_m}{ z-\gamma_m}\]  in
$\Hg(\Gamma, v).$ Then an appeal to the classical Riesz
representation theorem gives  that  for each sequence $(b_m) $ in $
\ell_v^2,$ we have
\begin{equation*}
 \sum_{m=1}^\infty a_m^n v_m \overline{b_m} \longrightarrow 0
\end{equation*} whenever $n\to \infty.$ Upon in particular taking
$b^{(l)}= \left(b_m^{(l)}\right)$ where $b_l^{(l)}= 1 $ and $0$
otherwise implies
\begin{equation}\label{limfact}
\lim_{n\to \infty}|a_m^n| v_m= 0
\end{equation} for each  $m.$ This consequence of the Riesz theorem will play
en essential role in the remaining part of the  proof.

We may first make the following splitting:
\begin{eqnarray}
\sum_{l=1}^\infty \int_{\Omega_l} |f_n(z)|^2d\mu(z)&\lesssim&
\sum_{l=1}^\infty
\int_{\Omega_l}\frac{1}{|z|^2}\left(\sum_{m=1}^{l-1}|a_m^n|v_m\right)^2
d\mu(z)\nonumber\\
&+& \sum_{l=1}^\infty \int_{\Omega_l}\frac{|a_l^n|^2
v_l^2}{|z-\gamma_l|^2}d\mu(z) + \sum_{l=1}^\infty \mu(\Omega_l)
\left(\sum_{m=l+1}^\infty\frac{|a_m^n|
v_m}{|\gamma_m|}\right)^2,\nonumber
\end{eqnarray} which  follows from the Cauchy--Schwarz inequality and
the growth condition \eqref{expon}.
 It suffices
to show that each of the three right-hand sums converges to zero
when $n\to \infty.$ We first show that
\begin{equation}
\label{triviallim} \lim_{n\to \infty} \sum_{l=1}^\infty
\int_{\Omega_l}\frac{|a_l^n|^2 v_l^2}{|z-\gamma_l|^2}d\mu(z) =0.
\end{equation}
From \eqref{trivialcomp}, for each small $\varepsilon>0,$ there
exists $N$ for which
\begin{equation*}
\int_{\Omega_l}\frac{ v_l}{|z-\gamma_l|^2}d\mu(z)< \varepsilon
\end{equation*} when $l> N.$ It follows that
\begin{eqnarray}
\label{commid} \sum_{l=1}^\infty\int_{\Omega_l}\frac{|a_l^n|^2
v_l^2}{|z-\gamma_l|^2}d\mu(z) &\lesssim& \sum_{l=1}^N |a_l^n|^2
v_l\int_{\Omega_l}\frac{ v_l}{|z-\gamma_l|^2}d\mu(z) +
\varepsilon\sum_{l=N+1}^\infty  |a_l^n|^2 v_l \nonumber\\
&\lesssim& \sum_{l=1}^N |a_l^n|^2 v_l +\varepsilon;
\end{eqnarray} here we used \eqref{nordeq} and \eqref{trivialcomp}. Taking the limit $n\to \infty$ in
\eqref{commid} and invoking \eqref{limfact} leads to the desired
conclusion \eqref{triviallim}.

 It  remains to prove
\begin{equation} \label{complimfinite}
\lim_{n\to \infty} \sum_{l=1}^\infty
\int_{\Omega_l}\frac{1}{|z|^2}\left(\sum_{m=1}^{l-1}|a_m^n|v_m\right)^2
d\mu(z)=0
\end{equation}
and
\begin{equation}
 \label{last}
 \lim_{n\to \infty} \sum_{l=1}^ \infty\mu(\Omega_l)
 \left(\sum_{m=l+1}^\infty\frac{|a_m^n| v_m}{|\gamma_m|}\right)^2=0.
 \end{equation}
For these limits, we only have to  modify the arguments  used to
establish
the estimates in  (3.2) and (3.3) 
in \cite{BMS2}.
 We set  \[\tau_l= \left( \int_{\Omega_l} |z|^{-2} d\mu(z)\right)^{\half}\] and first
 show that \eqref{complimfinite} holds. By duality, we
have
\begin{eqnarray*}
\left(\sum_{l=1}^\infty
\tau_l^2\left(\sum_{m=1}^{l-1}|a_m^n|v_m\right)^2 \right)^\half&=
&\sup_{\|c_l\|_\ell^2=1}
 \sum_{l=1}^\infty \tau_l |c_l|\sum_{m=1}^{l-1}|a_m^n|v_m\nonumber\\
& \leq & \sup_{\|c_l\|_\ell^2=1} \sum_{m=1}^\infty
|a_m^n|v_m\sum_{l=m+1}^{\infty} \tau_l |c_l|.\nonumber
\end{eqnarray*}
 Cauchy--Schwarz inequality applied to the last sum gives
\begin{equation}
\label{ineq} \left(\sum_{l=m+1}^{\infty} \tau_l |c_l|\right)^2 \leq
\sum_{l=m+1}^\infty \tau_l^2 V_l^\half\sum_{j= m+1}|c_j|^2
V_j^{-\half} .
\end{equation}
 By \eqref{compactcarleson}, we observe that
for each $\varepsilon >0,$ there exists $N_1$ for which $m\geq N_1,$
\begin{equation*}
\sum_{ l:2^k V_m < V_l\leq 2^{k+1} V_m} \tau_l^2 V_l^\half \lesssim
\frac{\varepsilon}{2^{k/2} V_{m+1}^\half}
\end{equation*} for $k\geq0$ and $m\geq N_1.$  Summing these inequalities for $m\geq N_1,$  we get
\begin{equation}
\label{taill} \sum_{l=m+1}^\infty \tau_l^2 V_l^\half \lesssim
\frac{\varepsilon}{V_{m+1}^\half}.
\end{equation}
Combining \eqref{ineq} with \eqref{taill}, we find
\begin{eqnarray}
\sum_{m=1}^\infty v_m \left(\sum_{l=m+1}^{\infty} \tau_l
|c_l|\right)^2&=& \sum_{m=1}^{N_1} v_m
 \left(\sum_{l=m+1}^{\infty} \tau_l |c_l|\right)^2+\sum_{m={N_1+1}}^\infty v_m
 \left(\sum_{l=m+1}^{\infty} \tau_l |c_l|\right)^2\nonumber\\
&\lesssim& \sum_{m=1}^{N_1} \frac{v_m}{V_{m+1}}\sum_{j= m+1}|c_j|^2
V_j^{-\half}
+ \varepsilon\underbrace{\sum_{m={N_1+1}}^\infty\frac{v_m}{V_{m+1}}
\sum_{j= m+1}^\infty|c_j|^2 V_j^{-\half}}_{\bigtriangleup}\nonumber\\
 &\lesssim& \sum_{m=1}^{N_1} \frac{v_m}{V_{m+1}}\sum_{j= m+1}|c_j|^2 V_j^{-\half}
  + \varepsilon.
  \label{cmpt}
  \end{eqnarray}
The  double sum in \eqref{cmpt} or $\bigtriangleup$ is bounded by an
absolute constant $C$ because
\begin{equation}
\sum_{m=1}^{N_1} \frac{v_m}{V_{m+1}}\sum_{j= m+1}|c_j|^2
V_j^{-\half}\leq \sum_{j= 1}^\infty |c_j|^2 V_j^{-\half}
\sum_{m=1}^{j-1} \frac{v_m}{V_{m+1}}
 \end{equation} when we change the order of the summation and so
\[ V_j^{-\half}\sum_{m=1}^{j-1} \frac{v_m}{V_{m+1}} \leq V_j^{-\half}
 \int_0^{V_j} t^{-\half} dt= 2. \]
 To obtain \eqref{complimfinite}, we see that
\begin{eqnarray*}
 \sum_{m=1}^\infty  |a_m^n|^2 v_m \sum_{m=1}^\infty\left(\sum_{l=m+1}^{\infty} \tau_l |c_l|\right)^2
 &\lesssim& C\sum_{m=1}^{N_1}  |a_m^n|^2 v_m
 + \varepsilon \sum_{m=N_1+1}^n  |a_m^n|^2 v_m \nonumber\\
 &\lesssim& \sum_{m=1}^{N_1}  |a_m^n|^2 v_m \longrightarrow 0
 \end{eqnarray*} as $n\to \infty $ which follows from \eqref{limfact}.

 To prove \eqref{last}, we note that the Cauchy--Schwarz inequality
gives
 \begin{equation*}
 \sum_{l=1}^\infty \mu(\Omega_l) \left(\sum_{m=l+1}^\infty\frac{|a_m^n| v_m}{|\gamma_m|}\right)^2 \leq
  \sum_{m=l+1}^\infty |a_m^n|^2
v_m P_{m-1}^{\half} \sum_{j=l+1}^\infty \frac{v_j}{P_{j-1}^{\half}
|\gamma_j|^2}.
\end{equation*}
 Since
\[ \sum_{j=l+1}^\infty
\frac{v_j}{P_{j-1}^{\half} |\gamma_j|^2}\le \int_0^{P_l}x^{-\half}
dx \le 2 P_l^{\half},\] it follows that
\[ \sum_{l=1}^\infty \mu(\Omega_l) \left(\sum_{m=l+1}^\infty
\frac{|a_m^n|v_m}{|\gamma_m|}\right)^2\lesssim \sum_{l=1}^\infty
\mu(\Omega_l) P_l^{\half} \sum_{m=l+1}^\infty |a_m^n|^2v_m
P_{m-1}^{\half},\] which becomes \[\sum_{l=1}^\infty \mu(\Omega_l)
\left(\sum_{m=l+1}^\infty
\frac{|a_m^n|v_m}{|\gamma_m|}\right)^2\lesssim \sum_{m=1}^\infty
|a_m^n|^2v_m P_{m-1}^{\half} \sum_{l=1}^{m-1} \mu(\Omega_l)
P_l^{\half}
\]
when we change the order of summation.
 By  \eqref{compactcarleson}, for each $\varepsilon >0, $
 there exists again an $N_2$
 for which for $m \geq N_2 $  it
follows that
\[\sum_{l: 2^k
P_{m-1}\leq P_l\leq 2^{k+1}P_{m-1}}\mu(\Omega_l) P_l^{\half}
 \lesssim
\frac{\varepsilon}{P^\half_{m-1} 2^{k/2}}.
\] Summing these inequalities  with respect to $k$ gives
\begin{equation*}
\sum_{l=1}^{m-1} \mu(\Omega_l) P_l^{\half} \lesssim
\frac{\varepsilon}{P^\half_{m-1} }
\end{equation*}
and we get
\begin{eqnarray*}
\sum_{m=1}^\infty |a_m^n|^2v_m P_{m-1}^{\half} \sum_{l=1}^{m-1}
\mu(\Omega_l) P_l^{\half} &\lesssim& \sum_{m=1}^{N_2} |a_m^n|^2v_m +
\varepsilon
\sum_{m=N_2+1}^\infty |a_m^n|^2v_m\nonumber\\
&\lesssim& \sum_{m=1}^{N_2} |a_m^n|^2v_m +
\varepsilon\longrightarrow 0
\end{eqnarray*} as $n\to \infty$ and completes the proof.

\subsection*{Acknowledgment}
The author would like to thank his PhD supervisor professor Kristian
Seip for his constructive  comments and encouragement during the
preparation of the paper.

\end{document}